\documentclass[draftcls,12pt, onecolumn]{IEEEtran}
\usepackage{amsfonts}
\usepackage{enumitem}
\usepackage{bbm}
\usepackage{times}
\usepackage{amssymb}
\usepackage{cleveref}
\usepackage{amsmath}
\usepackage{color,psfrag}
\usepackage{graphicx}
\usepackage{epsfig}
\usepackage{mathrsfs}
\usepackage{caption}
\usepackage{subfigure}
\usepackage{makeidx}
\usepackage{multirow}
\usepackage{dblfloatfix}
\usepackage{threeparttable}
\usepackage{dsfont}
\usepackage{slashbox}
\usepackage{lipsum}
\usepackage[font={footnotesize}]{caption}
\newlist{aims}{enumerate}{1}
\setlist[aims,1]{
  label={\arabic*.},
  leftmargin=0pt,
  align=left,
  labelsep=5pt,
  itemindent=\dimexpr\labelsep+\labelwidth+7pt\relax
}

\newtheorem{theorem}{Theorem}

\newtheorem{proposition}{Proposition}
\newtheorem{lemma}{Lemma}
\newtheorem{definition}{Definition}
\newtheorem{remark}{Remark}

\newtheorem{assumption}{Assumption}

\newcommand{\R}{\mathbb{R}}

\newcommand{\T}{\mathcal{T}}
\newcommand{\N}{\mathcal{N}}

\newcommand{\D}{\mathcal{D}}

\newcommand{\X}{\mathcal{X}}

\newcommand{\U}{\mathcal{U}}

\newcommand{\I}{\mathcal{I}}

\renewcommand{\L}{\mathcal{L}}

\renewcommand{\P}{\mathcal{P}}
\newcommand{\J}{\mathcal{J}}

\newcommand{\case}[1]{\begin{cases}#1\end{cases}}
\newcommand{\ar}[1]{\left[\begin{array}#1\end{array}\right]}
\newcommand{\al}[1]{\begin{align}#1\end{align}}
\newcommand{\eq}[1]{\begin{equation}#1\end{equation}}
\newcommand{\ald}[1]{\begin{aligned}#1\end{aligned}}

\begin{document}
\title{On Weak Topology for Optimal Control of Switched Nonlinear Systems}
\author{Hua Chen and Wei Zhang
\thanks{Hua Chen and Wei Zhang are with the Department of Electrical and Computer Engineering, The Ohio State University, Columbus, OH 43210. (e-mails: chen.3824@osu.edu; zhang.491@osu.edu)}}
\maketitle

\begin{abstract}
Optimal control of switched systems is challenging due to the discrete nature of the switching control input. The embedding-based approach addresses this challenge by solving a corresponding relaxed optimal control problem with only continuous inputs, and then projecting the relaxed solution back to obtain the optimal switching solution of the original problem. This paper presents a novel idea that views the embedding-based approach as a change of topology over the optimization space, resulting in a general procedure to construct a switched optimal control algorithm with guaranteed convergence to a local optimizer. Our result provides a unified topology based framework for the analysis and design of various embedding-based algorithms in solving the switched optimal control problem and includes many existing methods as special cases. 

\end{abstract}

%

\section{Introduction}
Switched systems consist of a family of subsystems and a switching signal determining the active subsystem (mode) at each time instant. Optimal control of switched systems involves finding both the continuous control input and the switching signal to jointly optimize certain system performance index. This problem has attracted considerable research attention due to its diverse engineering applications in power electronics~\cite{ONPDU09}, automotive systems~\cite{HR99,RDRK08,UDBPZ11}, robotics~\cite{WUZD13}, and manufacturing~\cite{CPW01}.  

Optimal control of switched systems is in general challenging due to the discrete nature of the switching control input, which prevents us from directly applying the classical optimal control techniques to solve the problem. To address this challenge, the maximum principle was extended in the literature to characterize optimal hybrid control solutions~\cite{P98,SC03,S99,S00}. However, it is still very difficult to numerically compute the optimal solutions based on these abstract necessary conditions~\cite{XA03}.  


Among the rich literature, one well-known method is the so-called bilevel optimization~\cite{XA03,XA04}. This approach divides the original optimal control problem into two optimization problems and solves them
at different levels. At the lower level, the approach fixes a switching mode sequence and optimizes the cost over the space of switching time instants through the classical variational approach.
At the upper level, the switching mode sequence is updated to optimize the cost. Although various heuristic schemes have been proposed for the upper level~\cite{EWA06,GVKSBT10,GVKSBT10_2}, solutions obtained via this method may still be unsatisfactory due to the restriction on possible mode sequences. 

More recently, an alternative approach based on the so-called \emph{embedding principle} has been proposed~\cite{BD05,VGBS13,VGBS13_2,WUZDB07}. This approach is closely related to the relaxed optimal control problems which optimize over the convex closure of the original control set. Several results concerning the existence property of the optimal solutions to the original problem have been discovered in the literature of relaxed optimal control problems~\cite{B74,GKNR96,W14}. 
%
The embedding-based approach adopts the idea of relaxing the control input and  takes one step further by introducing a projection operator which maps the relaxed optimal control back to the original input space to generate the desired switching control. There are three major steps involved in the embedding-based approach. The first step is to  embed the switched systems into a larger class of classical nonlinear systems with only continuous control inputs. Then, the optimal control of the relaxed system is obtained using the classical optimal control algorithms. Once the relaxed optimal solution is obtained, the solution to the original problem can be computed by projecting the relaxed solution back to the original input space through certain carefully designed projection operators. This approach has been successfully applied to numerous applications, such as power electronics~\cite{ONPDU09}, automotive systems~\cite{RDRK08,UDBPZ11}, and robotics~\cite{WUZD13}.

Several different versions of the embedding-based approach have been developed in the literature. These methods can be extended in their specific ways of embedding the switched trajectories, solving the associated classical optimal control problem, or projecting the relaxed solutions back to the original space. 
The main purpose of this paper is not trying to study these specific extensions by proposing different embedding-based algorithms. 
Instead, we aim to develop a general topology based framework for analyzing and designing various embedding-based optimal control algorithms.
The proposed framework is based on a novel idea that views the embedding-based approach as a change of topology over the optimization space. From this perspective, our framework adopts the weak topology structure and describes a general procedure to construct different switched optimal control algorithms.  
The framework involves several key components, and we derive conditions for these components under which the overall algorithm converges to a stationary point of the original switched optimal control problem. 
Our framework includes many existing results as special cases. We also illustrate the importance of viewing the switched optimal control problems from the topological perspective.

 The rest of this paper is organized as follows: Section~\ref{sec:probf} formulates the switched optimal control problems. 
 Section~\ref{sec:mainres} first reviews some important concepts in topology and then develops the proposed framework, along with its convergence analysis. 
 A numerical example demonstrating the use of our framework and the importance of selecting appropriate topology for particular problem is presented in Section~\ref{sec:exam}. Concluding remarks are given in Section~\ref{sec:concl}.



\section{Problem Formulation and Preliminaries}\label{sec:probf} 
Consider a switched nonlinear system model
\eq{\dot{x}(t) = f_{\sigma(t)}(t,x(t),u(t)), \text{ for a.e. } t\in [0,t_f],} where~$x(t)\in \R^{n_x}$ is the system state,~$u(t)\in U\subset \R^{n_u}$ is the continuous control input constrained in a compact and convex set $U$, and~$\sigma(t)\in \Sigma \triangleq \{1,2,\ldots,n_\sigma\}$ is the switching control input determining the active subsystem (mode) among a finite number $n_\sigma$ of subsystems at time~$t$. 

The cost function considered in our optimal control problem is given by~$h(x(t_f))$, i.e. only terminal state is penalized. Optimal control problems with nontrivial running costs can be transformed into this form by introducing an additional state variable~\cite{P97}.  It is also assumed that system~\eqref{eq:const} is subject to the following state constraints:
  \eq{h_j(x(t))\le 0, \quad \forall t\in [0,t_f], \quad \forall j\in \J \triangleq \{1,2,\ldots,n_c\}.\label{eq:const}}
The following assumption is imposed to ensure the existence and uniqueness of the state trajectory and the well-posedness of our optimal control problem.
\begin{assumption}\label{ass:lip}
\begin{aims}
\item~$f_i(t,x,u)$
is Lipschitz continuous with respect to all arguments for all~$i\in \Sigma$ with a common Lipschitz constant~$L$,
\item~$h_j(x)$,~$h(x)$
are Lipschitz continuous with respect to all argument for all~$j\in \J$ with a common Lipschitz constant~$L$.
\end{aims}
\end{assumption}
\begin{remark}
We assume a common Lipschitz constant $L$ to simplify notation. All the results in this paper extend immediately to the case where all these functions have different Lipschitz constants. 
\end{remark}

Following similar notations used in~\cite{BD05,VGBS13}, we rewrite the system dynamics as follows 
\eq{\ald{\dot{x} &= \sum\limits_{i=1}^{n_\sigma}d_i(t)f_i(t,x(t),u(t)) \triangleq f(t,x(t),u(t),d(t)), \text{ for a.e. }t\in [0,t_f],\label{eq:sysf}}}
where~$d(t) = [d_1(t),d_2(t),\ldots,d_{n_\sigma}(t)]\in D \triangleq \Big\{(d_1,\ldots,d_{n_\sigma}) \in \{0,1\}^{n_\sigma} \Big| \sum\limits_{i=1}^{n_\sigma} d_i = 1 \Big\}$ for a.e. $t\in[0,t_f]$, and $D$ is the set of corners of the~$n_\sigma$ simplex. The continuous input $u$ and switching input $d$ can be viewed as mappings from $[0,t_f]$ to $U$ and $D$, respectively. In this paper, we assume these mappings to be elements of the~$\L^2$ space, defined as follows.
\begin{definition}\label{def:l2space}
We say a function~$g:[0,t_f] \to G \subseteq \R^n$ belongs to~$\L^2([0,t_f],G)$, if \eq{\|g\|_{\L^2} \triangleq \left(\int\limits_{0}^{t_f} \|g(t)\|_2^2dt\right)^{\frac{1}{2}}<\infty,\label{eq:l2space}} where the integration is with respect to the Lebesgue measure. 
\end{definition}

 Let~$\U = \L^2([0,t_f],U)$ be the space of continuous control inputs and let~$\D = \L^2([0,t_f],D)$ be the space of switching control inputs. We denote by $\X_p = \U\times \D$ the overall original input space and call $\xi\in \X_p$ a \emph{original} input signal. Suppose the initial state $x(0) = x_0\in\R^{n_x}$ is given and fixed, we denote by~$x(t;\xi)$ the state at time $t$ driven by $\xi$ with initial state $x_0$. In order  to emphasize the dependence on~$\xi$, the following notations are adopted in this paper: \eq{\phi_t(\xi) \triangleq x(t;\xi), J(\xi ) \triangleq h(x(t_f;\xi)), \psi_{j,t}(\xi) \triangleq h_j(x(t;\xi)).\label{eq:newnote}}
We further define~$\Psi(\xi) \triangleq \max\limits_{j\in \J,t\in [0,t_f]}\big\{\psi_{j,t}(\xi)\big\}$ and the state constraints in \eqref{eq:const} can then be rewritten as $\Psi(\xi)\le 0$, since~$\Psi(\xi) = \max\limits_{j\in \J,t\in [0,t_f]}\big\{\psi_{j,t}(\xi)\big\} \le 0$ if and only if~$\psi_{j,t}(\xi)\le 0$ for all~$j\in \J$ and~$t\in [0,t_f]$.

Adopting the above notations, the optimal control problem of switched systems considered in this paper is reformulated as the following optimization problem:
\eq{\label{eq:pxp}\P_{\X_p}:  \case{&\inf\limits_{\xi\in \X_p} J(\xi),  \\ \text{subj. to }&\Psi(\xi) \le 0.} }

The problem $\P_{\X_p}$ is a constrained optimization problem over function space $\X_p$. However, the classical optimization techniques cannot be applied directly to solve this problem due to the discrete nature of~$\D$. The embedding-based approach is one of the most effective methods proposed in the literature for addressing this issue. This approach first embeds the switched systems into a larger class of traditional nonlinear systems with only continuous control inputs. Then, it solves an associated relaxed optimization problem through the classical numerical optimization algorithms. Lastly, it projects the relaxed optimal control back to the original input space to obtain the solution to the original problem. 
In this paper, we devise a novel idea that views the embedding-based approach as a change of topology of the optimization space, resulting in a general procedure for developing switched optimal control algorithms under the new topology.
In the next section, we first briefly review some concepts in weak topology and then establish the topology based framework.

\section{A Unified Framework for Switched Optimal Control Problem}\label{sec:mainres}
In this section, we establish the unified topology based framework to solve the switched optimal control problem $\P_{\X_p}$. We first introduce the weak topology to rigorously define the local minimizers of the problem. Then, the unified topology based framework is established and convergence of any algorithm constructed by the framework is proved provided that the conditions of the framework components are satisfied.

\subsection{Review of Weak Topology}\label{sec:mainresA}
Local minimizers are considered as solutions to general optimal control problems. Rigorous definition of local minimizers depends on the underlying topology adopted over the optimization space. 
Our framework adopts the weak topology which is defined as follows:

\begin{definition}[Weak Topology~\cite{R91}]
 Let~$\{g_i\}_{i\in \I}$ be a family of functions~$g_i: X\mapsto Y_i, \ \forall i\in \I$, mapping from a set~$X$ to several topological spaces~$Y_i$, respectively. The weak topology on~$X$ induced by~$\{g_i\}_{i\in \I}$, denoted by~$\T_{\{g_i\}_{i\in \I}}$, refers to the weakest topology on $X$ which makes all~$g_i$ continuous.
 \end{definition}
 \begin{remark}
 The structure of weak topology $\T_{\{g_i\}_{i\in \I}}$ is determined by the family of functions~$\{g_i\}_{i\in\I}$. In particular, the family may contain only one element. For example, the metric topology on a space $X$ is defined to be the weak topology induced by a norm function~$\|\cdot\|$, denoted by~$\T_{\|\cdot\|}$\footnote{Many norms can be defined on a space $X$ and each of them induces a metric topology. In this paper, we assume $\T_{\|\cdot\|_X} = \T_{\|\cdot\|_{\L^2}}$ if $X$ is a function space and $\T_{\|\cdot\|_X}=\T_{\|\cdot\|_2}$ if $X$ is an Euclidean space. }
 \end{remark}

The topology selected over the optimization space plays a critical role in characterizing local optimizers of the underlying optimization problem. Before providing the formal definition of a local minimizer, we first define a neighborhood around a point $\xi_p \in \X_p$ under a topology $\T_g$ as follows:
\begin{definition}\label{def:eneigh}
Given a topological space $(\X_p,\T_g)$, we say $\N_{\T_g}(\xi_p)\subset (\X_p,\T_g)$ is a neighborhood around $\xi_p$ under $\T_g$, if $\exists O\in \T_g$ such that $\xi_p\in O \subseteq \N_{\T_g}(\xi_p)$.
\end{definition}

Consider a mapping $g:\X_p\mapsto Y$, where $Y$ is a topological space endowed with a metric topology $\T_{\|\cdot\|_Y}$, a neighborhood around $\xi_p\in \X_p$ under $\T_g$ with radius $r$ is defined by:
 \eq{\N_{\T_g}(\xi_p,r) = \Big\{\xi_p'\in \X_p\Big| \|g(\xi_p) - g(\xi_p')\|_Y\le r\Big\}.}

Employing the above definition, a local minimizer of $\P_{\X_p}$ under a topology $\T_g$ is defined below. 
\begin{definition}\label{def:locminX}
We say $\xi_p^* \in \X_p$ is a local minimizer of~$\P_{\X_p}$ under the topology~$\T_g$, if there exists a neighborhood~$\N_{\T_g}(\xi_p^*)$ such that~$J(\xi_p^*)\le J(\xi_p'), \forall \xi_p'\in \N_{\T_g}(\xi_p^*)\cap \{\xi_p\in \X_p\big| \Psi(\xi_p)\le 0\}$.
\end{definition}

Different choices of the topologies $\T_g$ will lead to different characterization of local minimizers and hence affect the solution to the problem $\P_{\X_p}$. To further illustrate the importance of the weak topology in our framework, we provide a numerical example in Section~\ref{sec:exam} which shows that different topologies selected over the optimization space will result in different solutions to the same switched optimal control problem.

 \subsection{Solution Framework}\label{subsec:fc}

Note that it is difficult in general to directly check the local minimizer condition  in  Definition~\ref{def:locminX}, even for classical optimal control problems. In this paper, we adopt the following {\em optimality function} concept~\cite{P97} to characterize a necessary optimality condition. 
\begin{definition}\label{def:optfX}
A function $\theta_p(\cdot):\X_p\to \R$ satisfying the following conditions is called an optimality function for~$\P_{\X_p}$:
\begin{aims}
\item$\theta_p(\xi)\le 0$ for all~$\xi \in \X_p$;
\item if~$\xi_p^*$ is a local minimizer of~$\P_{\X_p}$, then~$\theta_p(\xi_p^*) = 0$.
\end{aims}
\end{definition}
\begin{remark}\label{rmk:theta}
Often times, the optimality function is required to be continuous (or upper semi-continuous)~\cite{P97}. Such a condition is introduced to ensure that in a topological space, if~$\xi^*$ is an accumulation point of any sequence~$\{\xi_i\}_{i\in \mathbb{N}}$ and~$\liminf\limits_{i\to \infty} \theta_p(\xi_i) = 0$, then we have~$\theta_p(\xi^*) = 0$. However, in our problem we do not assume the existence of accumulation points of the sequence~$\{\xi_i\}_{i\in \mathbb{N}}$. Hence, the continuity (or upper semi-continuity) condition is not necessary.
\end{remark}

Employing this optimality function definition and the necessary optimality condition encoded therein, our goal becomes constructing the optimization algorithm~$\Gamma_p:\X_p \to \X_p$ for $\P_{\X_p}$ such that~$\theta_p(\xi_p^i)\to 0$ as~$i\to \infty$, where $\{\xi_p^i\}_{i\in \mathbb{N}}$ is the sequence of original switched inputs generated by the optimization algorithm $\Gamma_p$ as defined in~\eqref{eq:iter} below.
\eq{\label{eq:iter}\xi_p^{i+1} = \case{ \Gamma_p(\xi_p^i), & \text{ if }\theta_p(\xi_p^i)<0,\\\xi_p^i, & \text{ if } \theta_p(\xi_p^i) = 0.}}
For simplicity, we denote by $\{\xi_p^i\}_{i\in\mathbb{N}}$ the sequence generated by~\eqref{eq:iter}.

Our topology-based framework involves three key steps and several important components given as follows. 
\begin{aims}
\item Relax the optimization space $\X_p$ to a vector space~$\X_r$, select a weak topology function~$g :\X_r\mapsto Y$ and construct a projection operator~$\mathscr{R}_k:\X_r\to\X_p$ associated with the weak topology~$\T_g$.
\item Solve the relaxed optimization problem~$\P_{\X_r}$ defined in~\eqref{eq:pxr} below by designing a relaxed optimality function~$\theta_r:\X_r\to\mathbb{R}$ and selecting (or constructing) a relaxed optimization algorithm~$\Gamma_r:\X_r\to\X_r$.
\item Set $\theta_p = \theta_r|_{\X_p}$ and~$\Gamma_p =\mathscr{R}_k\circ\Gamma_r$ with any initial condition $\xi^0_p\in \X_p$.

\end{aims}

The relaxed optimization problem~$\P_{\X_r}$ in the above framework is given by
\eq{\label{eq:pxr}\P_{\X_r}:  \case{&\inf\limits_{\xi\in \X_r} J(\xi), \\ \text{subj. to } &\Psi(\xi) \le 0,}}
and the relaxed optimality function $\theta_r$ is defined by replacing $\X_p$ and $\P_{\X_p}$ with $\X_r$ and $\P_{\X_r}$ in Definition~\ref{def:optfX}.

The main underlying idea of the proposed framework is to transform the switched optimization problem $\P_{\X_p}$ to a classical optimization problem $\P_{\X_r}$ which can be solved through the classical gradient-based methods in functional spaces~\cite{K70,P97}. The solution of $\P_{\X_r}$ will then be used to construct the solution to the original problem $\P_{\X_p}$. The key components of the framework include the relaxed optimization space $\X_r$, the weak topology $\T_g$, the projection operator $\mathscr{R}_k$, and the relaxed optimization algorithm characterized by $\theta_r$ and $\Gamma_r$. 

In the rest of this section, we will first show that $\theta_p$ is an optimality function for $\P_{\X_p}$ and then derive conditions for the aforementioned key components of our framework to guarantee that the sequence $\{\xi_p^i\}_{i\in\mathbb{N}}$ converges to a point satisfying the necessary optimality condition encoded in~$\theta_p$.

\subsection{Convergence Analysis and Proofs}

Before stating our main results, we first impose the following assumptions on $\X_r$, $\T_g$ and $\mathscr{R}_k$ in the framework to ensure its validity. 
\begin{assumption}\label{ass:stdass}
\begin{aims}
\item $J(\cdot)$ and~$\Psi(\cdot)$ are Lipschitz continuous under topology~$\T_g$ with a common Lipschitz constant~$L$.
\item $\X_p$ is dense in~$\X_r$ under~$\T_g$, i.e.~$\forall \xi_r\in \X_r$,~$\forall \epsilon>0$,~$\exists \xi_p\in \X_p$ s.t.~$\|g(\xi_r)-g(\xi_p)\|_Y \le \epsilon$.
\item There exists a projection operator~$\mathscr{R}_k:\X_r\to \X_p$ associated with~$\T_g$ and parametrized by~$k=1,2,\ldots$, such that~$\forall \xi_r\in\X_r$,~$\forall \epsilon>0$, there exists a~$\hat{k}\in \mathbb{N}$, such that  \eq{\left\|g(\mathscr{R}_k(\xi_r)) - g(\xi_r)\right\|_Y\le e_{\mathscr{R}_k}(k)\le \epsilon,\  \forall k\ge \hat{k}.\label{eq:projop}}
\end{aims}
\end{assumption}

Assumption~\ref{ass:stdass}.1 is a standard Lipschitz continuity condition that ensures the well-posedness of the relaxed problem $\P_{\X_r}$. Assumption~\ref{ass:stdass}.2 and~\ref{ass:stdass}.3 impose mild constraints on the relaxed space and topology that can be used in the framework. 

In the following lemma, we show that $\theta_p=\theta_r|_{\X_p}$ is an optimality function for $\P_{\X_p}$.
\begin{lemma}\label{lem:validtheta}
If $\theta_r$ is a valid optimality function for~$\P_{\X_r}$, then $\theta_p = \theta_r|_{\X_p}$ is a valid optimality function for~$\P_{\X_p}$.
\end{lemma}
\begin{IEEEproof}
To prove this lemma, we need to show~$\theta_p$ satisfies the conditions in Definition~\ref{def:optfX}. 
The first condition is trivially satisfied. For the second condition, suppose it does not hold, i.e. suppose~$\xi^*\in \X_p$  is a local minimizer for~$\P_{\X_p}$ but ~$\theta_p(\xi^*)< 0$. Since~$\theta_r(\xi^*)=\theta_p(\xi^*)$, by the definition of local minimizers for $\P_{\X_r}$, it follows that there exists a~$\xi_r$ and a positive number $C$, such that~$J(\xi_r)-J(\xi^*)\le -C$ and~$\Psi(\xi_r)\le -C$. By Assumption~\ref{ass:stdass}, we have~$\left|J(\mathscr{R}_k(\xi_r))-J(\xi_r)\right|\le L\left\| g(\mathscr{R}_k(\xi_r))-g(\xi_r)\right\|_Y\le L e_{\mathscr{R}_k}(k)$. By adding and subtracting $J(\xi_r)$, it follows that \eq{\ald{&J(\mathscr{R}_k(\xi_r))-J(\xi^*) \\ \le & \left|J(\mathscr{R}_k(\xi_r)) - J(\xi_r) \right| +J(\xi_r)-J(\xi^*)\\ \le & L e_{\mathscr{R}_k}(k)-C}}

For any given~$\xi_r\in \X_r$, choose~$\epsilon = \frac{C}{2L}$ in Assumption~\ref{ass:stdass}.3. For $k\ge \hat k $, it follows that~$L e_{\mathscr{R}_k}(k)-C\le-\frac{C}{2}<0$, hence~$J(\mathscr{R}_k(\xi_r))-J(\xi^*)<0$. A similar argument can be applied on~$\Psi$, yielding that~$\Psi(\mathscr{R}_k(\xi_r)) \le 0$. These statements contradict that~$\xi^*$ is a local minimizer for~$\P_{\X_p}$.
\end{IEEEproof}

%
%
To show the convergence of $\{\xi_p^i\}_{i\in \mathbb{N}}$, we adopt a similar idea of the sufficient descent property presented in~\cite{AWEV08}. In order to handle the projection step in our framework and the state constraints considered in our problem, we define two functions $Q:\X_r \times \mathbb{N} \mapsto \R$ and $P:\X_r \times \X_r\mapsto \R$ below.

\eq{\label{eq:G} \ald{Q(\xi,k) \triangleq \case{\ald{\max\{&J(\mathscr{R}_k\circ \Gamma_r(\xi)) - J(\Gamma_r(\xi)),\\ &\Psi(\mathscr{R}_k \circ \Gamma_r (\xi))-\Psi(\Gamma_r(\xi))\},} &  \text{if } \Psi(\xi)\le 0,\\  \Psi(\mathscr{R}_k \circ \Gamma_r(\xi)) - \Psi(\Gamma_r(\xi)),  &  \text{if } \Psi(\xi)>0.}}}
\eq{P(\xi_1,\xi_2) \triangleq \case{\max\{J(\xi_2)-J(\xi_1), \Psi(\xi_2)\}, &  \text{if } \Psi(\xi_1)\le 0,\\  \Psi(\xi_2)-\Psi(\xi_1), &  \text{if } \Psi(\xi_1)>0,}}

We introduce the function $Q$ to compactly characterize the change of  the cost $J$ and the constraint $\Psi$ at a point $\xi$ under the projection operator $\mathscr{R}_k$. For a feasible point, we care about both the changes of the cost and the constraint under $\mathscr{R}_k$. For an infeasible point, we only care about the change of the constraint.

The function $P$ characterizes the value difference of $J$ and $\Psi$ between two points $\xi_1$ and $\xi_2$. If $\xi_1$ is feasible and $P<0$, it means the cost can be reduced while maintaining feasibility by moving from $\xi_1$ to $\xi_2$. Similarly, if $\xi_1$ is infeasible and $P<0$, it is possible to reduce the infeasibility by moving from $\xi_1$ to $\xi_2$.

Exploiting Assumption~\ref{ass:stdass}.3, a bound for the function $Q$ is derived in the following lemma.
\begin{lemma}\label{lem:rk}
There exists a~$k^*\in\mathbb{N}$ such that given~$\omega\in (0,1)$, for any~$C>0$,~$\gamma_C>0$, and for any~$\xi\in\X_p$ with~$\theta_p(\xi)<-C$, we have  \eq{Q(\xi,k)\le (\omega-1)\gamma_C\theta_p(\xi) ,\   \forall k\ge k^*.}
\end{lemma}
\begin{IEEEproof}
This is a straightforward result from Assumption~\ref{ass:stdass}.1,  Assumption~\ref{ass:stdass}.2 and Lemma~\ref{lem:validtheta}.
\end{IEEEproof}

Employing the definition of the function $P$ and the above two lemmas, our main result on the convergence of $\{\xi_p^i\}_{i\in \mathbb{N}}$ is presented below.
\begin{theorem}\label{thm:mainres}
If for each~$C>0$, there exists a~$\gamma_C>0$ such that for any~$\xi_r\in \X_r$ with~$\theta_r(\xi_r)<-C$, 
\eq{P(\xi_r,\Gamma_r(\xi_r))\le \gamma_C\theta_r(\xi_r)<0,} 
then for an appropriate choice of $k$ for $\mathscr{R}_k$, for any $\xi_p^0\in \X_p$ the following two conclusions hold:
\begin{aims}
\item if there exists a~$i_0\in \mathbb{N}$ such that~$\Psi(\xi_p^{i_0})\le 0$, then~$\Psi(\xi_p^i)\le 0$ for all~$i\ge i_0$,
\item$\lim\limits_{i\to \infty}\theta_p(\xi_p^i) = 0$, i.e. the sequence $\{\xi_p^i\}_{i\in\mathbb{N}}$ converges asymptotically to a stationary point.
\end{aims}
\end{theorem}
 
\begin{IEEEproof}
\begin{aims}
\item Suppose there exists an~$i_0$ such that~$\Psi(\xi_p^{i_0})\le 0$, then we have for $k\ge k^*$ \eq{\ald{&\Psi(\xi_p^{i_0+1}) \\  = &  \Psi(\mathscr{R}_k(\Gamma_r(\xi_p^{i_0}))) - \Psi(\Gamma_r(\xi_p^{i_0}))+\Psi(\Gamma_r(\xi_p^{i_0}))-  \Psi(\xi_p^{i_0}) +\Psi(\xi_p^{i_0}) \\ \le& (\omega-1)\gamma_C\theta_p(\xi_p^{i_0}) +\gamma_C\theta_r(\xi_p^{i_0})\\ =& \omega\gamma_C\theta_p(\xi_p^{i_0})<0}}
\item We need to consider two cases due to different form of~$P$ for different values of~$\Psi$.
\begin{itemize}
\item Case 1:~$\Psi(\xi_p^i)>0$ for all~$i\in \mathbb{N}$, i.e. the entire sequence is infeasible. 

Suppose~$\lim\limits_{i\to \infty}\theta_p(\xi_p^i)\neq 0$, since~$\theta_p(\cdot)$ is a non-positive function, we know there must exists~$C>0$ such that~$\liminf\limits_{i\to \infty}\theta_p(\xi_p^i) = -2C$. Hence, there exists an infinite subsequence~$\{\xi_p^{i_m}\}$ and an~$m_1 \in \mathbb{N}_+$ such that~$\theta_p(\xi_p^{i_m})<-C$ for all~$m\ge m_1$. Then, it follows that for all~$m\ge m_1$, and for $k\ge k^*$, we have
 \eq{\ald{&\Psi(\xi_p^{i_{m+1}})- \Psi(\xi_p^{i_m})\\  =& \Psi(\mathscr{R}_k\circ\Gamma_r(\xi_p^{i_m})) -\Psi(\Gamma_r(\xi_p^{i_m}))+\Psi(\Gamma_r(\xi_p^{i_m}))- \Psi(\xi_p^{i_m})\\  \le& (\omega-1)\gamma_C\theta_p(\xi_p^{i_m}) +\gamma_C\theta_r(\xi_p^{i_m})\\ =& \omega\gamma_C\theta_p(\xi_p^{i_m})<0}}
This leads to the fact that~$\liminf_{m\to \infty} \Psi(\xi_p^{i_m}) = -\infty$, which contradicts the lower boundedness of~$\Psi$  implied by Assumption~\ref{ass:lip}.
\item Case 2: There exists an~$i_0$ such that~$\Psi(\xi_p^{i_0})\le 0$. 

By the first conclusion, it follows that~$\Psi(\xi_p^i)\le 0$ for all~$i\ge i_0$. 
 Suppose~$\liminf\limits_{i\to \infty}\theta_p(\xi_p^i)\neq 0$, then there exists~$C>0$  such that~$\liminf\limits_{i\to \infty}\theta_p(\xi_p^i) = -2C$. Hence, there exists an infinite subsequence~$\{\xi_p^{i_m}\}$ and a ~$m_1 \in \mathbb{N}_+$ such that~$\theta_p(\xi_p^{i_m})<-C$ for all~$m\ge m_1$. Then, it follows that for all~$m\ge m_1$ and for all~$k\ge k^*$, we have:
\eq{\ald{&J(\xi_p^{i_{m+1}})- J(\xi_p^{i_m})\\=& J(\mathscr{R}_k\circ \Gamma_r^l(\xi_p^{i_m}))- J(\xi_p^{i_m})\\ =& J(\mathscr{R}_k\circ \Gamma_r^l(\xi_p^{i_m})) - J(\Gamma_r^l(\xi_p^{i_m}))+J(\Gamma_r^l(\xi_p^{i_m})) - J(\xi_p^{i_m}) \\ \le& (\omega-1)\gamma_C\theta_p(\xi_p^{i_m}) +\gamma_C\theta_r(\xi_p^{i_m})\\ = &\omega\gamma_C\theta_p(\xi_p^{i_m})<0 }}
This leads to the fact that~$\liminf\limits_{m\to \infty} J(\xi_p^{i_m}) = -\infty$, which contradicts with the lower boundedness of~$J$ implied by Assumption~\ref{ass:lip}.
\end{itemize}
\end{aims}
\end{IEEEproof}

In the following section, a concrete numerical example is shown to illustrate the use of the proposed framework and the importance of viewing the switched optimal control problem from the topological perspective. 

\section{Illustrating Example}\label{sec:exam}
 Numerous embedding-based switched optimal control algorithms proposed in the literature can be analyzed using the proposed framework. Depending on the underlying applications, one may choose different relaxed spaces $\X_r$, weak topologies $T_g$, optimality functions $\theta_r$, projection operators $\mathscr{R}_k$, or relaxed optimization algorithms $\Gamma_r$. Each combination of these components will lead to a different switched optimal control algorithm that may have a better performance for particular problems. 

In this section, we present a numerical example to illustrate how the proposed framework can be used to guide the design and analysis of a switched optimal control algorithm.
In addition, we will also show through the example that proper selection of the weak topology is important for obtaining a satisfactory solution. 

\begin{figure}[ht] 
\centering
\includegraphics[width=0.4\linewidth]{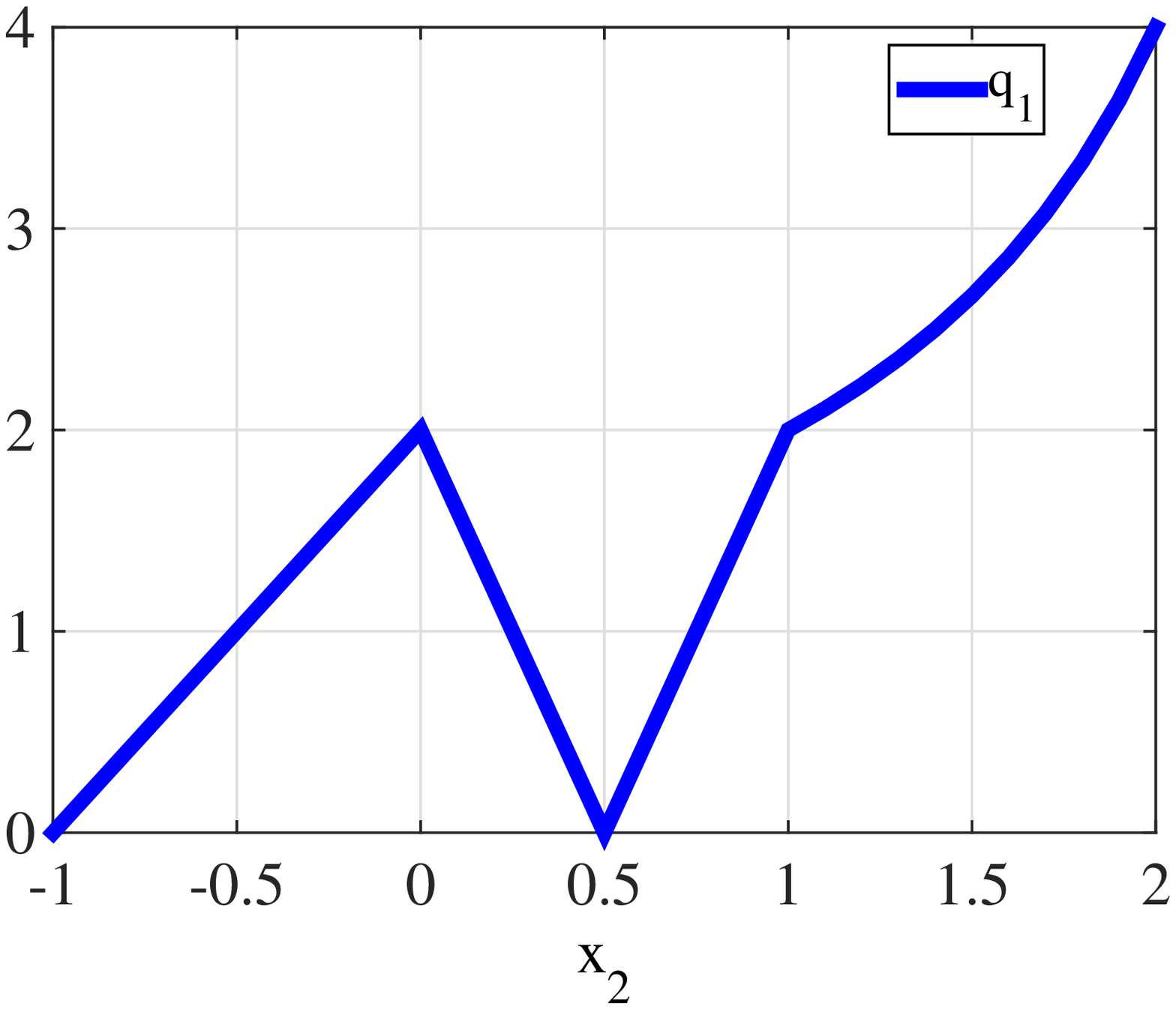}
\includegraphics[width=0.4\linewidth]{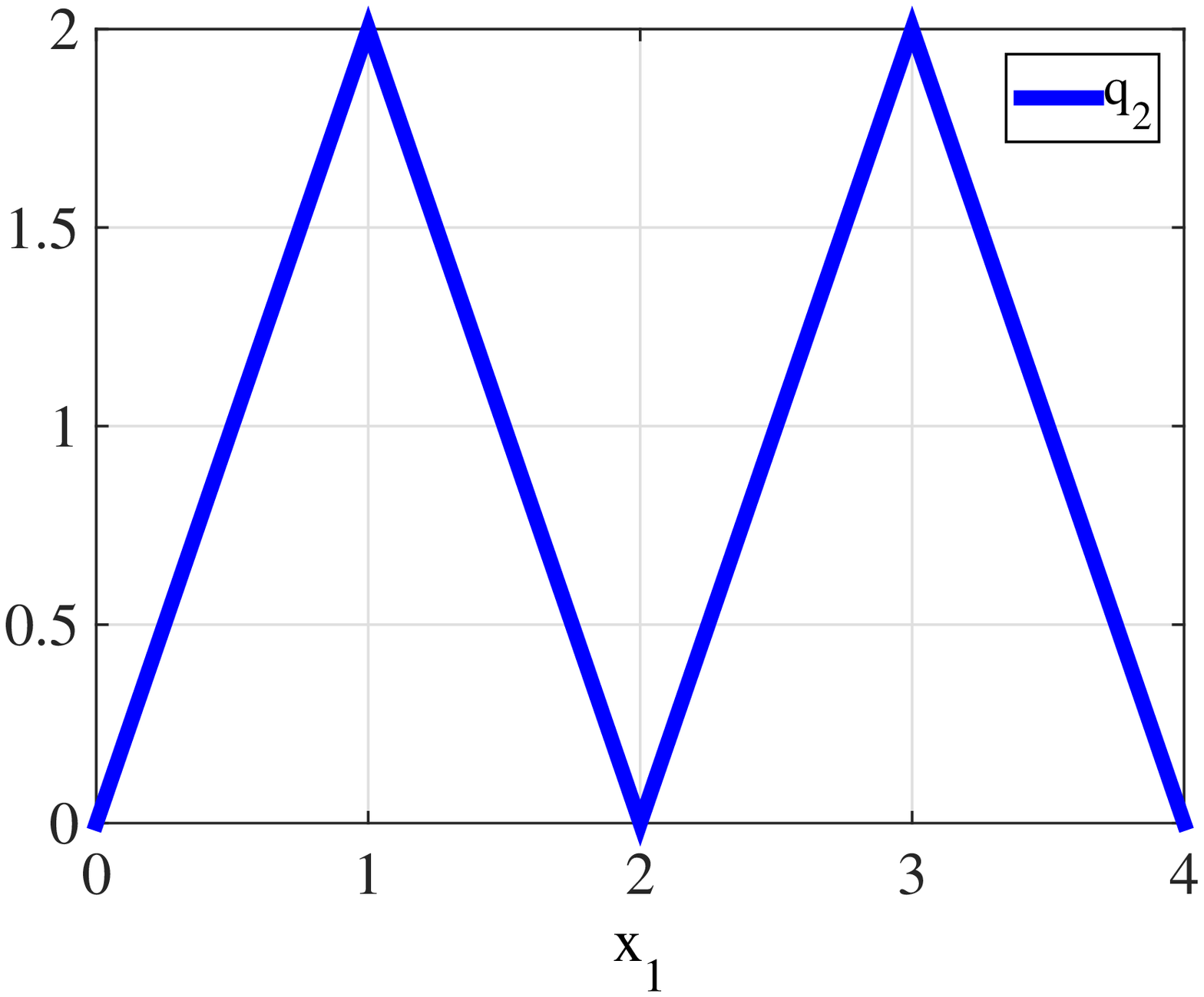} 
\caption{ Illustration of the functions $q_1(x_2)$ and $q_2(x_1)$}
\label{fig:vecg}
\end{figure}

Consider the following switched system consisting two subsystems in the domain given by~$ \Big\{x=(x_1,x_2)^T\in \R^2\Big\}$. Dynamics of each mode is given by: 
\eq{\ald{\text{ Mode 1: }\dot x &= f_1(x_1,x_2) = \ar{{cc}q_1(x_2), &0}^T,  \\ \text{ Mode 2: }\dot x &= f_2(x_1,x_2) = \ar{{cc}0, &q_2(x_1)}^T,}} where~$q_1$ and~$q_2$ are defined by~\eqref{eq:gfun} as follows and are illustrated in Fig.~\ref{fig:vecg}. Suppose, for simplicity, that neither continuous input nor state constraints are involved in our problem and denote the control signal by $\xi(t) = (d_1(t),d_2(t))^T$ where $d_1(t)$ and $d_2(t)$ are the discrete inputs defined in~\eqref{eq:sysf}. Let $x(0) = [0,0]^T$ be the initial state and let the time horizon be~$[0,2]$. The cost function is given by~$h(x(2;\xi)) = \|x(2;\xi) - A\|_2$ where~$A= [3,2]^T$. In other words, we want to find the optimal switching input $\xi$ to  minimize the Euclidean distance between the terminal state and point $A$. It is not difficult to see that any input signal resulting in terminal state at $A$ is a global minimizer of this problem with the optimal cost of $0$.


\al{q_1(x_2)& = \case{ 0 , & \text{ if } x_2\le 0, \\ 2x_2+2, & \text{ if } x_2 \in [-1,0),\\-4x_2+2 & \text{ if } x_2\in [0,0.5),\\4x_2-2 & \text{ if } x_2\in [0.5,1),\\ \frac{4}{3-x_2}, & \text{ if } x_2\in[1,2],\\ 4, & \text{ if } x_2>2.} &  q_2(x_1) &= \case{0 , & \text{ if } x_1\le 0, \\2x_1, & \text{ if } x_1\in[0,1),\\-2x_1+4, & \text{ if } x_1\in[1,2),\\ 2x_1-4, & \text{ if } x_1\in [2,3),\\-2x_1+8, &\text{ if } x_1\in [3,4],\\ 0, & \text{ if } x_1>4.}\label{eq:gfun}}

To utilize our framework, we first reformulate the optimal control problem as an optimization problem over function space described in~\eqref{eq:pxp} as follows.
\eq{\min\limits_{\xi\in \X_p} J(\xi),\label{eq:simuopt}} where $\X_p = \L^2([0,2],D)$ is the optimization space with $D = \{(d_1,d_2)\in \{0,1\}^2\Big| d_1+d_2=1\}$ and~$J(\xi) = \|\phi_{2}(\xi) - A\|_2$ with $\phi_2(\xi)$ adopting the notation introduced by~\eqref{eq:newnote}. 

We can apply the existing algorithm developed in~\cite{VGBS13,VGBS13_2} to solve this switched optimal control problem where the weak topology is chosen to be the one induced by the entire state trajectory. For this numerical example, the algorithm is terminated whenever the optimality function is sufficiently close to zero. The detailed termination condition is given by $\theta(\xi_p^k)>-\epsilon$, where
$\epsilon$ is chosen to be $10^{-6}$. We discretize the time horizon $[0,2]$ into $N = 2^8 = 64$ samples as $\{t_i = \frac{i-1}{N}\}_{i=1}^{N}$ and let the initial state $\xi_p^0 = (d_1^0,d_2^0)^T$ be the switching input signal defined by \al{\label{eq:ini1}d_1^0(t) &= \case{1, & \text{ if } t\in [0,t_{50}], \\ 0, & \text{ if } t\in(t_{50},2].}  \\ \label{eq:ini2}d_2^0(t) &= 1-d_1^0(t_i), \forall t\in[0,2].}
\begin{figure}[htbp] 
\centering
\includegraphics[width=0.6\linewidth]{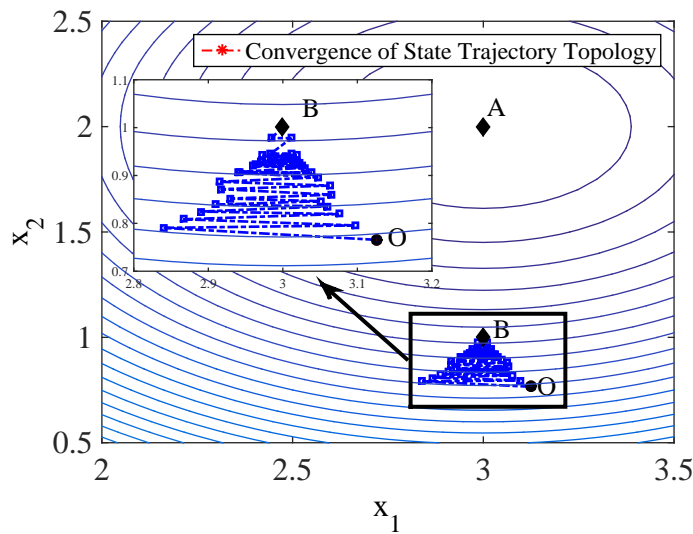}
\caption{ Convergence of the terminal states generated by the switched optimal control algorithm under the weak topology induced by the entire state trajectory}
\label{fig:simu1}
\end{figure}

Fig.~\ref{fig:simu1} shows the convergence of the terminal states of the trajectories generated by the algorithm developed in~\cite{VGBS13}. In the figure, the black solid circle $O$ is the terminal state generated by the initial input signal $\xi^0$ defined in~\eqref{eq:ini1} and~\eqref{eq:ini2}. Point $A$ is the terminal state corresponding to the global minimizer which is shown as a black diamond. It is clear that the solution obtained through this algorithm converges to a stationary point with terminal state at $B = [3,1]^T$ which is also shown as a black diamond. The cost associated with the solution is $1$ which does not equal the cost of the global minimizer of the problem. This is because the neighborhood of any local minimizer under the weak topology induced by the state trajectory excludes those switching inputs which generate close enough terminal states but not close enough entire state trajectories.

Since in this particular problem, the cost depends only on the terminal state and no constraints are inovlved, it is nature to consider the weak topology induced by the terminal state. In the following, we will use such a weak topology and propose a modified projection operator and other components in the framework. The proposed framework can directly be used to analyze the convergence of the new algorithm. We now detail the choices as follows.

\begin{aims}
\item $\X_r$: $\X_r = \L^2([0,2],D_r)$, where~$D_r \triangleq \Big\{(d_1,d_2) \in [0,1]^2 \Big| d_1+d_2= 1 \Big\}$, which adopts the general idea of taking the convex closure of the original input space~$D$.
\item $\T_g$:  $\T_g$ is chosen to be  the weak topology induced by the terminal state function~$\phi_2(\cdot)$ and is denoted by~$\T_{\phi_2}$.
\item $\mathscr{R}_k$: $\mathscr{R}_k$ is the frequency modulation operator with frequency~$2^k$ which is defined as follows. 
Let $\xi^R = \mathscr{R}_k(\xi) = (d^R_1,d^R_2)^T\in \X_p$ be the projected signal given by
\eq{\label{eq:proj_def} \xi^R(t) = \case{1, & \text{ if } t\in (T_{i,1}, T_{i,2}), \\ 0, & \text{ otherwise,} }\quad  \forall i = 1,2,\ldots,k,}
where $T_{i,j}$ is given by:
\eq{\label{eq:Timescale} \ald{T_{i,1} &= \frac{i-1}{2^{k-1}}+ \frac{1}{2}\int\limits_{t_i}^{t_{i+1}}1-d_1(t)dt= \frac{i-1}{2^{k-1}}+ \frac{1}{2}\int\limits_{t_i}^{t_{i+1}}d_2(t)dt,\\ T_{i,2} &= T_{i,1}+\int\limits_{t_i}^{t_{i+1}}d_1(t)dt, \\ T_{i,3} &= \frac{i}{2^{k-1}},}}
where $\{t_i\}_{i=1}^{k} = \{\frac{i}{2^{k-1}}\}_{i=1}^{k}$ is a partition of the time horizon $[0,2]$. 
\item $\theta_r$: $\theta_r(\xi) = \min\limits_{\xi_r\in \X_r}DJ(\xi,\xi_r-\xi)$, where~$DH(x;x') = \lim\limits_{\lambda\downarrow 0}\frac{H(x+\lambda x')-H(x)}{\lambda}$ is the directional derivative for function $H$ at $x$ along direction $x'$.

\item $\Gamma_r$: $\Gamma_r$ is chosen to be the gradient descent optimization algorithm given by: $\Gamma_r = \hat{\Gamma}^l$ where $\hat{\Gamma}$ is the standard steepest decent algorithm and $l$ is determined by verifying the condition of Theorem~\ref{thm:mainres}, i.e. for any $\xi\in \X_r$, $l$ is determined as follows:
\eq{ l = \min\{k\in \mathbb{N} \Big| \ J(\hat{\Gamma}^l(\xi)) - J(\xi)\le \gamma_C\theta_r(\xi)\},} where $\gamma_C$ is the constant in Theorem~\ref{thm:mainres}. 

\end{aims}

\begin{proposition}
The components specified as above satisfy the conditions of the topology based framework, i.e. given any initial condition, the sequence of switched inputs generated by the algorithm $\Gamma_p = \mathscr{R}_k\circ \Gamma_r$ converges to a stationary point of this problem.
\end{proposition}
\begin{IEEEproof}
To prove this proposition, it only needs to be shown that $\theta_r$ is a valid optimality function of the relaxed problem and  Assumption~\ref{ass:stdass} and the condition in Theorem~\ref{thm:mainres} are satisfied by the above choices. 
\begin{itemize}
\item Validity of $\theta_r$: 
\begin{itemize}
\item $\theta_r(\xi) =  \min\limits_{\xi_r\in \X_r}DJ(\xi,\xi_r-\xi) \le DJ(\xi,\xi - \xi) = 0$; 
\item Suppose $\xi$ is a local minimizer of $\P_{\X_r}$ but $\theta_r(\xi)<0 $, then $\exists \xi'$ such that $DJ(\xi;\xi'-\xi)<0$. By mean value theorem, we have $\exists \lambda\in(0,1)$ such that $J(\xi+\lambda(\xi'-\xi))-J(\xi)  = \lambda DJ(\xi;\xi'-\xi)+ o(\lambda) < 0$. 
\end{itemize}
\item Assumption 2.1:  For any two switched input $\xi^1$ and $\xi^2$, we have
\eq{\ald{\left|J(\xi^1) - J(\xi^2)\right| &= \left|\|\phi_{2}(\xi^1) - A\|_2 - \|\phi_{2}(\xi^2) - A\|_2 \right| \\ &\le \|\phi_{2}(\xi^1) - \phi_{2}(\xi^2)\|_2},}
where the last inequality is due to the triangle inequality and the Lipschitz constant can be taken to be $1$. Since there is no constraint in this problem, Assumption 2.1 is satisfied. 
\item Assumption 2.2:  By the chattering lemma~\cite{BD05, B74}, $\X_p$ is dense in $\X_r$ under the weak topology induced by  the entire state trajectory $\T_\phi$ which is stronger than $\T_{\phi_2}$. Hence $\X_p$ is dense in $\X_r$ under the weaker topology $\T_{\phi_2}$ induced by the terminal state.
\item Assumption 2.3:  The validity of this projection operator is ensured by an analogous argument of the proof of Theorem 1 in~\cite{BD05}. 

\item Condition in Theorem~\ref{thm:mainres}: This is clearly satisfied due to our construction of $\Gamma_r$.
\end{itemize}
\end{IEEEproof}



%

\begin{figure}[htbp] 
\centering
\includegraphics[width=0.6\linewidth]{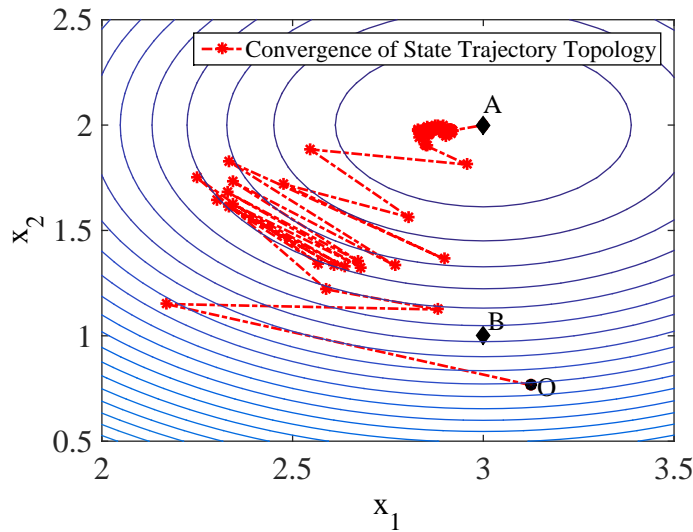}
\caption{ Convergence of the terminal states generated by the switched optimal control algorithm under the weak topology induced by the terminal state}
\label{fig:simu2}
\end{figure}
We implement the algorithm developed as above to solve the switched optimal control problem with the same initial settings as previous. Fig.~\ref{fig:simu2} shows the convergence of the terminal states of the trajectories generated by this algorithm. The resulting sequence generated by this algorithm actually converges to the global minimizer with terminal state at $A = [3,2]^T$ and the associated cost is given by $0$.  This is because under the new weak topology (induced by the terminal state), the solution obtained through the algorithm in~\cite{VGBS13} is not a stationary point anymore.  Therefore, the weak topology induced by the terminal state is more appropriate than the weak topology induced by the entire state trajectory for this particular problem.

This numerical example shows how our framework can be used for analyzing and designing various switched optimal control algorithms and the importance of choosing appropriate weak topology for different underlying problems.

\section{Conclusion}\label{sec:concl}
In this paper, we present a unified topology based framework that can be used for designing and analyzing various embedding-based switched optimal control algorithms.

Our framework is based on a novel viewpoint which considers the embedding-based methods as a change of topology over the optimization space. From this viewpoint, our framework adopts the weak topology structure and develops a general procedure to construct a switched optimal control algorithm. Convergence property of the algorithm is guaranteed by specifications on several key components involved in the framework. A concrete numerical example is provided to demonstrate the use of the proposed framework and the importance of selecting the appropriate weak topology in our framework. 

Possible extensions of this work include the considering the switched optimal control problems with switching costs and other forms of constraints. 


\bibliographystyle{IEEEtranS}
\bibliography{HChen}
\end{document}